\documentclass{tsp-short}

\newtheorem{lem}{Lemma}[section]

\theoremstyle{remark}
\newtheorem{rem}{Remark}[section]
\theoremstyle{definition}
\newtheorem{dfn}{Definition}[section]

\begin{document}

\title
[System of interacting particles with Markovian switching] {System
of interacting particles with Markovian switching}

\author{A. A. Pogorui}
\address{A. A. Pogorui, Zhitomir State Franko University, Velyka Berdychivska 40, Zhitomir 10008, Ukraine}
\email{pogor@zu.edu.ua}


\date{10/04/2013}

\begin{abstract}
Most of the published articles on random motions have been devoted
to the study of the telegraph process or its generalizations that
describe the random motion in $R^n$ of a single particle in a
Markov or semi-Markov medium. However, up to our best knowledge
there are no published papers dealing with the interaction of two
or more particles which move according to the telegraph processes.
In this paper, we construct the system of telegraph processes with
interactions, which can be interpreted as a model of ideal gas. In
this model, we investigate the free path times of a family of
particles before they are collided with any other particle. We
also study the distribution of particles, which described by
telegraph processes with hard collisions and reflecting
boundaries, and investigate its limiting properties.
\end{abstract}

\maketitle

\section{Introduction}

Let $\left\{\xi \left(t\right),\ t\ge 0\right\}$ be a Markov
process on the phase space $\left\{0,1\right\}$ with generative
matrix
\[Q=\lambda \left( \begin{array}{cc}
-1 &  1 \\
 1 & -1 \end{array} \right).\]

\begin{dfn} $S\left(t\right)$ is the telegraph process if
\begin{equation*}
\frac{d}{dt}S\left(t\right)=v{\left(-1\right)}^{\xi
\left(t\right)},\ \ v=const>0,
\end{equation*}
\[S\left(0\right)=y_0.\]
\end{dfn}
For a set of real numbers $y_1<y_2<\dots <y_n$ we consider a
family of independent telegraph processes $S_i\left(t\right)$,
$i=1,2,\dots,n$ with $S_i\left(0\right)=y_i$. It is assumed that
all the processes have absolute velocity $v$ and parameter of
switching process $\lambda>0$, and starting form $y_i$ the process
$S_i(t)$ has equal probabilities of initial directions of the
motion.

Denote by $x(y_i,t)$ the position of the particle $i$ at time $t$,
which starts from site $y_i$. Suppose that particle $x(y_i,t)$
develops as the telegraph process $S_i(t)$ up to the hard
collision with another particle. Under the hard collision of two
particles, we mean that at the time of the collision the particles
change their direction to the opposite that is, the particles
exchange the telegraph processes that describe their movement. It
is easily verified that the positions of the particles $x(y_i,t)$,
$i=1,2,\dots,n$ at time $t$ coincide with the order statistics of
$S_i (t)$, $i=1,2,\dots,n$ as follows

\begin{equation}\label{1}
x(y_1,t)=S_{(1)}(t), x(y_2,t)=S_{(2)}(t),\dots, x(y_n,t)=S_{(n)}
(t). \end{equation}

\begin{rem}
It should be noted that each $x(y_i,t)$, $i=1,2,\dots,n$ is not a
telegraph process for all $t\geq 0$.
\end{rem}

\begin{rem}It follows from the description of
$x\left(y_i,t\right)$ that $x\left(y_1,t\right)\le
x\left(y_2,t\right)\le \cdots \le x\left(y_n,t\right)$ for any
$t\geq 0$. Such kind of model for Wiener processes with
coalescence after collision are called the Arratia flow and they
were studied in \cite{p6}-\cite{p8}.
\end{rem}

\noindent Various problems such as the number of particle
collisions up to time $t$ in the Arratia flow are studied in
\cite{p9}.

Below the explicit form for the distribution of the meeting
instant of two telegraph processes on the line, which started at
the same time from different positions in the line, is obtained.
We also study the limiting distribution of the meeting instant of
two telegraph processes on the line under Kac's condition. It
allows us to investigate the system of telegraph processes with
interactions, which can be interpreted as a model of ideal gas. In
this model, we investigate the free path times of a family of
particles before they are collided with any other particle. We
also study the distribution of particles, which described by
telegraph processes with hard collisions and reflecting
boundaries, and investigate its limiting properties.

\section{Distribution of the first collision of two telegraph particles}

Consider two particles $1$ and $2$ on a line. Each particle can
move in two opposite directions. Starting at $x_i\in R$, $i=1,2$
particle $i$ moves at the velocity $v>0$ in one of two directions
during a random time interval that is exponential distributed with
parameter $\lambda >0$. Then the particle changes its direction
and so on. In the sequel, such particle is said to be a telegraph
particle as its motion satisfies the telegraph equation \cite{p1},
\cite{p2}.

\noindent Let ${\xi }_1\left(t\right)$, ${\xi }_2\left(t\right)$
be independent alternating Markov processes with the phase space
$\left\{0,1\right\}$ and with the generator matrix $Q$.

\noindent Denote by $x_i(t)$ the position of particle $i$ at a
point of time $t\geq0$ up to the first collision with another
particle. It is easily seen that
\[\frac{d}{dt}x_i\left(t\right)=v{\left(-1\right)}^{{\xi }_i\left(t\right)},\]
\[x_i\left(0\right)=x_i.\]
We assume that $z=x_2-x_1>0$ and put
$\Delta\left(t\right)=x_2\left(t\right)-x_1\left(t\right)$.

\noindent Denote $\eta \left(t\right)=\left({\xi
}_1\left(t\right),{\xi }_2\left(t\right)\right)$. Suppose $\eta
\left(0\right)=\left(k_1,k_2\right)$ and define
\[{\tau }_{\left(k_1,k_2\right)}\left(z\right)=inf\left\{t\geq 0:\ \Delta\left(t\right)=0\right\},\ \ k_j\in \left\{0,1\right\}.\]

\noindent Denote by
$f_{\left(k_1,k_2\right)}\left(t,z\right)dt=P\left({\tau
}_{\left(k_1,k_2\right)}\left(z\right)\in dt\right)$ the density
probability function (pdf) of ${\tau
}_{\left(k_1,k_2\right)}\left(z\right)$.

\begin{lem} \label{l1}
For $t\geq\frac{z}{2v}$

\begin{equation}\label{2}
f_{\left(0,1\right)}\left(t,z\right)=e^{-2\lambda t}\delta
\left(z-2vt\right)+\frac{z\lambda }{2v^2}e^{-2\lambda
t}\frac{I_1\left(\frac{\lambda
}{v}\sqrt{4v^2t^2-z^2}\right)}{\sqrt{4v^2t^2-z^2}},\
\end{equation}

\begin{equation}\label{3}
f_{\left(0,0\right)}\left(t,z\right)=f_{\left(1,1\right)}\left(t,z\right)=\frac{z\lambda
}{2v^2}e^{-2\lambda
t}\int^t_{{z}/{2v}}{\frac{I_1\left(\frac{\lambda
}{v}\sqrt{4v^2u^2-z^2}\right)}{\sqrt{4v^2u^2-z^2}}\frac{I_1\left(2\lambda
\left(t-u\right)\right)}{\left(t-u\right)}}du
\end{equation}

\noindent and

\begin{equation}\label{4}
f_{\left(1,0\right)}\left(t,z\right)=\frac{z\lambda
}{2v^2}e^{-2\lambda
t}\int^t_{{z}/{2v}}{\frac{I_1\left(\frac{\lambda
}{v}\sqrt{4v^2u^2-z^2}\right)}{\sqrt{4v^2u^2-z^2}}\int^{t-u}_0{\frac{I_1\left(2\lambda
\left(t-u-v\right)\right)}{\left(t-u-v\right)}\frac{I_1\left(2\lambda
v\right)}{v}dv}}du.
\end{equation}
\end{lem}

\noindent \emph{Proof.} Let us consider the Laplace transforms of
${\tau }_{\left(k_1,k_2\right)}\left(z\right)$, $k_i\in
\left\{0,1\right\}$.
\[{\varphi }_{\left(k_1,k_2\right)}\left(s,z\right)=E\left[e^{-s{\tau }_{\left(k_1,k_2\right)}\left(z\right)}\right],\ \ s>0.\]

\noindent By using the renewal theory, we can obtain the following
system of integral equations for these Laplace transforms:

\begin{align*}
{\varphi }_{\left(0,1\right)}\left(s,z\right) = &
e^{-\frac{s+2\lambda }{2v}z}+\frac{\lambda
}{2v}\int^z_0{e^{-\frac{s+2\lambda }{2v}u}{\varphi
}_{\left(1,1\right)}\left(s,z-u\right)du}\\ +& \frac{\lambda
}{2v}\int^z_0{e^{-\frac{s+2\lambda }{2v}u}{\varphi
}_{\left(0,0\right)}\left(s,z-u\right)du}\\ = &
e^{-\frac{s+2\lambda }{2v}z}+\frac{\lambda
}{2v}e^{-\frac{s+2\lambda }{2v}z}\int^z_0{e^{\frac{s+2\lambda
}{2v}u}\left({\varphi
}_{\left(1,1\right)}\left(s,u\right)+{\varphi
}_{\left(0,0\right)}\left(s,u\right)\right)du}
\end{align*}
\noindent
\begin{align*}
{\varphi }_{\left(0,0\right)}\left(s,z\right) = & \lambda
\int^\infty_0{e^{-\left(s+2\lambda \right)u}{\varphi
}_{\left(0,1\right)}\left(s,z\right)du\ }+\lambda
\int^\infty_0{e^{-\left(s+2\lambda \right)u}{\varphi
}_{\left(1,0\right)}\left(s,z\right)du\ } \\
 = & \frac{\lambda }{s+2\lambda }\left({\varphi
 }_{\left(0,1\right)}\left(s,z\right)+{\varphi
}_{\left(1,0\right)}\left(s,z\right)\right),
\end{align*}
\noindent
\begin{align*}
{\varphi }_{\left(1,1\right)}\left(s,z\right) = & \lambda
\int^\infty_0{e^{-\left(s+2\lambda \right)u}{\varphi
}_{\left(0,1\right)}\left(s,z\right)du\ }+\lambda
\int^\infty_0{e^{-\left(s+2\lambda \right)u}{\varphi
}_{\left(1,0\right)}\left(s,z\right)du\ } \\ = & \frac{\lambda
}{s+2\lambda }\left({\varphi
}_{\left(0,1\right)}\left(s,z\right)+{\varphi
}_{\left(1,0\right)}\left(s,z\right)\right),
\end{align*}
\noindent
\begin{align*}
{\varphi }_{\left(1,0\right)}\left(s,z\right) = & \frac{\lambda
}{2v}\int^\infty_0{e^{-\frac{s+2\lambda }{2v}u}\left({\varphi
}_{\left(0,0\right)}\left(s,z+u\right)+{\varphi
}_{\left(1,1\right)}\left(s,z+u\right)\right)du\ } \\ = &
\frac{\lambda }{2v}e^{\frac{s+2\lambda
}{2v}z}\int^\infty_z{e^{-\frac{s+2\lambda }{2v}u}\left({\varphi
}_{\left(0,0\right)}\left(s,u\right)+{\varphi
}_{\left(1,1\right)}\left(s,u\right)\right)du}.
\end{align*}
It is easily seen that

\begin{equation}\label{5}
{\varphi }_{\left(0,0\right)}\left(s,z\right)={\varphi
}_{\left(1,1\right)}\left(s,z\right).
\end{equation}

\noindent Taking into account that

\begin{equation}\label{6}
{\varphi }_{\left(0,0\right)}\left(s,z\right)+{\varphi
}_{\left(1,1\right)}\left(s,z\right)=\frac{2\lambda }{s+2\lambda
}\left({\varphi }_{\left(0,1\right)}\left(s,z\right)+{\varphi
}_{\left(1,0\right)}\left(s,z\right)\right),
\end{equation}

\noindent we have
\[\frac{\partial }{\partial z}{\varphi }_{\left(0,1\right)}\left(s,z\right)=-\frac{s+2\lambda }{2v}{\varphi }_{\left(0,1\right)}\left(s,z\right)+\frac{{\lambda }^2}{v\left(s+2\lambda \right)}\left({\varphi }_{\left(0,1\right)}\left(s,z\right)+{\varphi }_{\left(1,0\right)}\left(s,z\right)\right),\]
\[\frac{\partial }{\partial z}{\varphi }_{\left(1,0\right)}\left(s,z\right)=\frac{s+2\lambda }{2v}{\varphi }_{\left(1,0\right)}\left(s,z\right)-\frac{{\lambda }^2}{v\left(s+2\lambda \right)}\left({\varphi }_{\left(0,1\right)}\left(s,z\right)+{\varphi }_{\left(1,0\right)}\left(s,z\right)\right).\]

\noindent It is well-known \cite{p3} that ${\varphi
}_{\left(1,0\right)}\left(s,z\right)$ and ${\varphi
}_{\left(0,1\right)}\left(s,z\right)$ satisfy the following
equation
\[det\left( \begin{array}{cc}
\frac{\partial }{\partial z}+\frac{\left(s+2\lambda \right)}{2v}-\frac{{\lambda }^2}{v\left(s+2\lambda \right)} & -\frac{{\lambda }^2}{v\left(s+2\lambda \right)} \\
\frac{{\lambda }^2}{v\left(s+2\lambda \right)} & \frac{\partial
}{\partial z}-\frac{\left(s+2\lambda \right)}{2v}+\frac{{\lambda
}^2}{v\left(s+2\lambda \right)} \end{array}
\right)f\left(z\right)=0.\]

\noindent By calculating the determinant, we get
\[\frac{{\partial }^2}{\partial z^2}f\left(z\right)-\frac{s^2+4\lambda s}{4v^2}f\left(z\right)=0.\]
Solving this equation, we have
\[f\left(z\right)=C_1e^{\sqrt{s^2+4\lambda s}\ \frac{z}{2v}}+C_2e^{-\sqrt{s^2+4\lambda s}\ \frac{z}{2v}}.\]

\noindent The constants obtained from the system of integral
equations yields
\begin{equation}\label{7}
{\varphi}_{\left(0,1\right)}\left(s,z\right)=e^{-\frac{z}{2v}\sqrt{s^2+4\lambda
s}}
\end{equation}
and
\begin{equation}\label{8}
{\varphi }_{\left(1,0\right)}\left(s,z\right)=\frac{s+2\lambda
-\sqrt{s^2+4\lambda s}}{s+2\lambda +\sqrt{s^2+4\lambda
s}}e^{-\frac{z}{2v}\sqrt{s^2+4\lambda s}}.
\end{equation}

\noindent Taking into account Eqs.(\ref{5}),(\ref{6}), we have
\begin{equation}\label{9}
{\varphi }_{\left(0,0\right)}\left(s,z\right)={\varphi
}_{\left(1,1\right)}\left(s,z\right)=\frac{s+2\lambda
-\sqrt{s^2+4\lambda s}}{2{\lambda
}}e^{-\frac{z}{2v}\sqrt{s^2+4\lambda s}}.
\end{equation}
The inverse Laplace of ${\varphi
}_{\left(0,1\right)}\left(s,z\right)$ yields the following pdf
(\cite{p4}, p.239, Formula 88)

\begin{align*}
f_{\left(0,1\right)}\left(t,z\right) = & {{\mathcal
L}}^{-1}\left(e^{-\frac{z}{2v}\sqrt{s^2+4\lambda s}},t\right)\\ =
& e^{-2\lambda t}\delta \left(z-2vt\right)+2z\lambda e^{-2\lambda
t}\frac{I_1\left(\frac{\lambda
}{v}\sqrt{4v^2t^2-z^2}\right)}{\sqrt{4v^2t^2-z^2}},\ \ \ t\geq
\frac{z}{2v}.
\end{align*}

\noindent Hence, Eq.(\ref{2}) is proved and
\[P\left({\tau }_{\left(0,1\right)}\left(z\right)\in dt\right)=e^{-2\lambda t}\delta \left(z-2vt\right)dt+2z\lambda e^{-2\lambda t}\frac{I_1\left(\frac{\lambda }{v}\sqrt{4v^2t^2-z^2}\right)}{\sqrt{4v^2t^2-z^2}}dt.\]

\noindent It is easily verified that

\begin{equation}\label{10}
\exp{\left\{-\frac{z}{2v}\sqrt{s^2+4\lambda
s}\right\}}=\exp{\left\{-\frac{z}{2v}s+\int^{\infty
}_0{\left(1-e^{-sy}\right)\frac{\lambda }{v}\frac{{e^{-2\lambda
y}I}_1\left(2\lambda y\right)}{y}dy}\right\}}.
\end{equation}

\noindent Then, it comes from \cite{p12}, p.237, no.49 the
following inverse Laplace transform

\begin{equation*}
{{\mathcal L}}^{-1}\left(1+\frac{s-\sqrt{s^2+4\lambda s}}{2\lambda
},t\right)=e^{-2\lambda t}\frac{I_1\left(2\lambda t\right)}{t}.
\end{equation*}

\noindent It is easily seen that the following condition holds

\begin{equation}\label{11}
\int^{\infty }_0{\left(1\wedge t\right)}e^{-2\lambda
t}\frac{I_1\left(2\lambda t\right)}{t}dt<+\infty .
\end{equation}

It is well known that the distribution which the Laplace transform
can be represented as the right side of Eq.(\ref{10}) with the
condition (\ref{11}) belongs to the infinitely divisible
distribution \cite{p13}. Therefore, the pdf
$f_{\left(0,1\right)}\left(t,z\right)$ is the infinitely divisible
density function.

\noindent Using  \cite{p12}, p.237, no.49, we get
\begin{align*}{{\mathcal L}}^{-1}\left(\frac{s+2\lambda -\sqrt{s^2+4\lambda s}}{s+2\lambda +\sqrt{s^2+4\lambda s}},t\right)= & \frac{1}{{4\lambda }^2}\ {{\mathcal L}}^{-1}\left({\left(s+2\lambda -\sqrt{s^2+4\lambda s}\right)}^2,t\right)\\
= & e^{-2\lambda t}\int^t_0{\frac{I_1\left(2\lambda
\left(t-v\right)\right)}{\left(t-v\right)}\frac{I_1\left(2\lambda
v\right)}{v}dv}\end{align*}
By calculating
\begin{align*}
{{\mathcal L}}^{-1}\left(\frac{s+2\lambda
-\sqrt{s^2+4\lambda s}}{s+2\lambda +\sqrt{s^2+4\lambda
s}}e^{-\frac{z}{2v}\sqrt{s^2+4\lambda s}}\right),
\end{align*}

\noindent we obtain Eq.(\ref{4}).

It is easily seen that $f_{\left(0,1\right)}\left(t,z\right)$ is a
heavy tail probability density function w.r.t. $t$. Indeed, by
using asymptotic expansion for $I_1\left(t\right)$ \cite{p5}, we
have

\begin{equation}\label{12}
{\mathop{\lim }_{t\to +\infty }\sqrt{2\pi t} I_1\left(t\right)\
}e^{-t}=1.
\end{equation}

\noindent Therefore,
\[E{\left[{\tau }_{\left(0,1\right)}\left(z\right)\right]}^{\alpha }\ge 2z\lambda \int^\infty_{\frac{z}{2v}}{t^{\alpha }e^{-2\lambda t}\frac{I_1\left(\frac{\lambda }{v}\sqrt{4v^2t^2-z^2}\right)}{\sqrt{4v^2t^2-z^2}}dt}=+\infty,\ \ {\rm for}\ \alpha \ge \frac{1}{2}.\]
It is easily verified that $E{\left[{\tau
}_{\left(0,1\right)}\left(z\right)\right]}^{\alpha }<\infty $ for
$0\le \alpha <\frac{1}{2}$.

For ${\tau }_{(0,1)}(z)$ at time $t=0$ particles move in opposite
directions to meet each other and for ${\tau }_{(0,1)}(z)$ at time
$t=0$ particles move in opposite directions far away from each
other.

\noindent Hence, $E{\tau }_{\left(0,1\right)}\left(z\right)\le \
E{\tau }_{\left(1,0\right)}\left(z\right)$ and
$f_{\left(1,0\right)}\left(t,z\right)$ is also a heavy tail
density function w.r.t. $t$.

\noindent Let us consider the following so-called Kac's condition
(or the hydrodynamic limit): denote by $\lambda ={\varepsilon
}^{-2}$, $v=c{\varepsilon }^{-1}$, as $\varepsilon >0$, that is
$v\rightarrow+\infty$, and $\lambda \rightarrow +\infty$, such
that $\frac{v^2}{\lambda }\rightarrow c^2$.

\noindent It was proved in \cite{p1} that under Kac's condition
the telegraph process $x\left(t\right)$ weakly converges to the
Wiener process $w\left(t\right)\sim N\left(0,c^2t\right).$

\noindent Denote $f\left(t,z\right)=\frac{cz\ {\rm
exp}\left(-\frac{c^2z^2}{4t}\right)}{2\sqrt{\pi }t^{{3}/{2}}}.$ It
is well known that $f(t,z)$ is the pdf of a collision instant of
two particles moving according to Wiener paths $w(t)$, where $z>0$
is the distance between starting points of the particles.

\begin{lem} \label{l2}
For each $k_1,k_2\in \{0,1\}$, $f_{(k_1,k_2)}(t,z)$ weakly
converges to $f(t,z)$ under Kac's condition.
\end{lem}

\noindent \emph{Proof.} It follows from Eqs.(\ref{7})-(\ref{9})
that
\[{\mathop{\lim }_{\varepsilon \rightarrow 0} {\varphi }_{\left(k_1,k_2\right)}\left(s,z\right)\ }=e^{-zc\sqrt{s}}.\]
Passing to the inverse Laplace transform, we have
\[f\left(t,z\right)={{\mathcal L}}^{-1}\left(e^{-zc\sqrt{s}}\right)=\frac{cz\ {\rm exp}\left(-\frac{c^2z^2}{4t}\right)}{2\sqrt{\pi }t^{{3}/{2}}}.\]

Therefore, under Kac's conditions not only the telegraph process
weakly converges to the Wiener process, but the first meeting
instant of two telegraph processes weakly converges to the first
meeting instant of the corresponding two Wiener processes.

\begin{rem}
It should be noted that instead of two telegraph processes
$x(y_1,t)$, $x(y_2,t)$ on the line we can consider the bivariate
process $\overrightarrow{x}(t)=(x(y_1,t),x(y_2,t))$ on the plane.
The process $\overrightarrow{x}(t)$ is driven by the switching
process $\eta(t)$. Denote $l=\{(x,y): x=y; x,y \in \mathbb{R}\}$.
For this case
\[{\tau }_{\left(k_1,k_2\right)}\left(z\right)=inf\left\{t\geq 0:\ \overrightarrow{x}(t)\in l\right\}\]
\end{rem}

\section{Estimation of the number of particle collisions.}

\noindent Denote by $N_{\left(0,1\right)}\left(t,z\right)$ the
number of collisions of particles $x_i\left(t\right),\ i=1,2$
during time $(0,t)$, $t>0$ assuming $\eta
\left(0\right)=\left(0,1\right)$.

\noindent Consider the renewal function
$H_{\left(0,1\right)}\left(t,z\right)=EN_{\left(0,1\right)}\left(t,z\right)$.
By using the Laplace transform for general renewal function
\cite{p16}, it follows from Eqs.(\ref{7})-(\ref{8}) that the
Laplace transform
${\hat{H}}_{\left(0,1\right)}\left(s,z\right)={\mathcal
L}\left(H_{\left(0,1\right)}\left(t,z\right),s\right)$ of
$H_{\left(0,1\right)}\left(t,z\right)$ w.r.t. $t$ has the
following form

\begin{align*}
{\hat{H}}_{\left(0,1\right)}\left(s,z\right)= &
\frac{e^{-\frac{z}{2v}\sqrt{s^2+4\lambda s}}}{s}\sum^{\infty
}_{k=0}{{\left(\frac{s+2\lambda -\sqrt{s^2+4\lambda s}}{s+2\lambda
+\sqrt{s^2+4\lambda s}}\right)}^k} \\ = &
e^{-\frac{z}{2v}\sqrt{s^2+4\lambda s}}\left(\frac{s+2\lambda
+\sqrt{s^2+4\lambda s}}{2s\sqrt{s^2+4\lambda s}}\right).
\end{align*}

\noindent It is easily verified that
\[{{\mathcal L}}^{-1}\left(\frac{s+2\lambda +\sqrt{s^2+4\lambda s}}{2s\sqrt{s^2+4\lambda s}}\right)=\frac{1}{2}+\left(\left(\frac{1}{2}+\lambda t\right)I_0\left(2\lambda t\right)+\lambda tI_1\left(2\lambda t\right)\right)e^{-2\lambda t}.\]

\noindent Therefore,

\begin{align}\label{13}
\nonumber H_{\left(0,1\right)}\left(t\right)= &
\int^t_{\frac{z}{2v}}{e^{-2\lambda u}\left(\delta
\left(z-2vu\right)+2z\lambda \frac{I_1\left(\frac{\lambda
}{v}\sqrt{4v^2u^2-z^2}\right)}{\sqrt{4v^2u^2-z^2}}\right)} \\
\nonumber \times & \left(\frac{1}{2}+e^{-2\lambda
\left(t-u\right)} \left(\left(\frac{1}{2}+\lambda
\left(t-u\right)\right)I_0\left(2\lambda
\left(t-u\right)\right)+\lambda \left(t-u\right)I_1\left(2\lambda
\left(t-u\right)\right)\right)\right)du \\ = & \nonumber
\frac{e^{-\frac{\lambda z}{v}}}{2}+z\lambda
\int^t_{\frac{z}{2v}}{e^{-2\lambda u}\frac{I_1\left(\frac{\lambda
}{v}\sqrt{4v^2u^2-z^2}\right)}{\sqrt{4v^2u^2-z^2}}}du
\\
+ & e^{-2\lambda t}\left(\left(\frac{1}{2}+\lambda
\left(t-\frac{z}{2v}\right)\right)I_0\left(2\lambda
\left(t-\frac{z}{2v}\right)\right)+\lambda
\left(t-\frac{z}{2v}\right)I_1\left(2\lambda
\left(t-\frac{z}{2v}\right)\right)\right)\\
\nonumber + & e^{-2\lambda t}z\lambda
\int^t_{\frac{z}{2v}}{\frac{I_1\left(\frac{\lambda
}{v}\sqrt{4v^2u^2-z^2}\right)}{\sqrt{4v^2u^2-z^2}}}
\\
\nonumber \times & \left(1+\left(\left(1+2\lambda
\left(t-u\right)\right)I_0\left(2\lambda
\left(t-u\right)\right)+\lambda \left(t-u\right)I_1\left(2\lambda
\left(t-u\right)\right)\right)\right)du.
\end{align}
\\
It follows from Eq.(\ref{13}) that by putting $\lambda
={\varepsilon }^{-2}$, $v=c{\varepsilon }^{-1}$, we have
\[H_{\left(0,1\right)}\left(t,z\right)=O\left({\varepsilon }^{-1}\right)=O\left(\sqrt{\lambda }\right)=O\left(v\right)\ \ {\rm as}\ \varepsilon \to 0.\]
For ${y,y}^*$ such as ${y<y}^*$ and a fixed $T>0$ denote by
$\widetilde{\tau} =inf\left\{T;
t:x\left(y,t\right)-x\left(y^*,t\right)=0\right\}$.

\noindent Much in the same way we can show that for all
$k_1,k_2\in \left\{0,1\right\}$

\[H_{\left(k_1,k_2\right)}\left(t,z\right)=O\left({\varepsilon }^{-1}\right)=O\left(\sqrt{\lambda }\right)=O\left(v\right)\ \ {\rm as}\ \varepsilon \to 0.\]

\begin{lem} \label{l3}
There exist $C>0$ such that for any two points $y,y^*$ ($y<y^*$)
\[E\widetilde{\tau }\leq C\left(y^*-y\right).\]
\end{lem}

\noindent \emph{Proof.}
\begin{align*}
E{\tau }_{\left(0,1\right)}= & \int^T_{\frac{z}{2v}}{te^{-2\lambda
t}\left[\delta \left(z-2vt\right)+2z\lambda
\frac{I_1\left(\frac{\lambda
}{v}\sqrt{4v^2t^2-z^2}\right)}{\sqrt{4v^2t^2-z^2}}\right]dt}
\\ \leq & \frac{z}{2v}+2z\lambda
\int^T_{\frac{z}{2v}}{t\frac{I_1\left(\frac{\lambda
}{v}\sqrt{4v^2t^2-z^2}\right)}{\sqrt{4v^2t^2-z^2}}dt}\\ = &
\frac{z}{2v}+\frac{z}{2v}\left(I_0\left(\frac{\lambda
}{v}\sqrt{4T^2v^2-z^2}\right)-1\right)\leq Cz,
\end{align*}

\noindent where $C=\frac{I_0\left(2T\lambda \right)}{2v}$.

\noindent Now
\begin{align*}
E{\tau }_{\left(1,0\right)}= &
\int^T_{{z}/{2v}}{tf_{\left(1,0\right)}\left(t,z\right)dt}\\
= & 2z\lambda \int^T_{{z}/{2v}}{te^{-2\lambda
t}\int^t_{{z}/{2v}}{\frac{I_1\left(\frac{\lambda
}{v}\sqrt{4v^2u^2-z^2}\right)}{\sqrt{4v^2u^2-z^2}}\int^{t-u}_0{\frac{I_1\left(2\lambda
\left(t-u-r\right)\right)}{\left(t-u-r\right)}\frac{I_1\left(2\lambda
r\right)}{r}dr}}dudt}\\ < &Cz, \end{align*}

\noindent where $C=\frac{\lambda }{v}\int^T_0{te^{-2\lambda
t}\int^t_0{\frac{I_1\left(2\lambda
u\right)}{u}\int^{t-u}_0{\frac{I_1\left(2\lambda
\left(t-u-r\right)\right)}{\left(t-u-r\right)}\frac{I_1\left(2\lambda
r\right)}{r}dr}}dudt}$.

\section{Free path times of a family of particles.}

Since we consider the model of an ideal gas it is natural to
assume that the number of particles is very large. As an example,
we consider a model with an infinite number of particles, and
study the free path of the particles before collisions.

\noindent Consider the segment $\left[0,S\right]\subset
\mathbb{R}$ and an increasing sequence of different points
$\left\{y_n;\ n\geq 1\right\}$ from this segment. As above, we
consider a family of independent telegraph processes $S_k(t)$,
$k\geq1$ and trajectories $x(y_k,t)$ of particles, which satisfy
Eq.(\ref{1}).

Introduce the following random times:
\[{\tau }_1=T>0,\]
\[{\tau }_k=inf\left\{T;t:{\left(x\left(y_k,t\right)-x\left(y_{k-1},t\right)\right)}=0\right\},\ \ k\geq 2.\]
The random variable ${\tau }_k$ is a time of free path of the
particle with number $k$ up to the collision with a particle
starting with a smaller number or until $T$ (finite) if no one
collision occurs.

\begin{lem} \label{l4}
Suppose $\left\{y_n;\ n\geq 1\right\}\subset \left[0,S\right]$,
$0<S<+\infty $ is a sequence of different points. Then
\[\sum^\infty_{k=1}{{\tau }_k}<+\infty \ \ {\rm a}.{\rm s}.\]
\end{lem}

\noindent \emph{Proof.} Consider the following random times
\[{\widetilde{\tau}}_1=T,\]
\[{\widetilde{\tau}}_k=inf\left\{T;t:{\left(S_k\left(t\right)-S_{k-1}\left(t\right)\right)}=0\right\},\ \ k\geq 2.\]
It is easily seen that ${\tau }_k\leq {\widetilde{\tau}}_k$ for
all $k\geq1$.

Hence, if we show that that
$\sum^\infty_{k=1}{{\widetilde{\tau}}_k}<+\infty \ \ {\rm a}.{\rm
s}.$, we prove the lemma. Since ${\widetilde{\tau}}_k\ge 0$, it is
sufficient to prove that
\[\sum^{\infty }_{k=1}{E{\widetilde{\tau}}_k}<+\infty .\ \ \]
Consider the set of numbers $y_1<y_2<\dots <y_n$. It follows from
Lemma \ref{l1} that there exists $C>0$ such that for any $k\geq 2$
\[E{\widetilde{\tau}}_k\leq C\left(y_{k}-y_{k-1}\right).\]
Hence, we have

\begin{align}\label{14}
{\mathop{\lim }_{n \rightarrow \infty}
\sum^n_{k=2}{E{\widetilde{\tau}}_k}\ }\leq C{\mathop{\lim }_{n
\rightarrow \infty} \sum^n_{k=2}{\left(y_{k}-y_{k-1}\right)}\
}=C\sum^\infty_{k=2}{\left(y_{k}-y_{k-1}\right)}\leq CS.
\end{align}

Therefore, it follows from Eq.(\ref{14}) that
$\sum^\infty_{k=1}{{\widetilde{\tau}}_k}$ converges almost surely
and it concludes the proof.

\noindent Note that Lemma \ref{l4} for Wiener particles was proved
in \cite{p8}.

Let us denote by $N_{\left(k_1,k_2,\dots
,k_n\right)}\left(t,y_1{\rm ,\ }y_2,\dots ,y_n\right),\ \ k_i\in
\left\{0,1\right\}$, $y_1<y_2<\dots <y_n$ the number of collisions
of particles $x\left(y_i,t\right),\ i=1,2,\dots ,n$ during time
$(0,t)$, $t>0$ assuming $\eta \left(0\right)=\left(k_1,k_2,\dots
,k_n\right)$.

\noindent Then it is easily seen that
\begin{eqnarray*}
H_{\left(k_1,k_2,\dots ,k_n\right)}\left(t,y_1{\rm ,\ }y_2,\dots
,y_n\right)&=& EN_{\left(k_1,k_2,\dots ,k_n\right)}\left(t,y_1{\rm
,\ }y_2,\dots
,y_n\right)\\
&=&\sum_{i=1}^{n-1}{H_{\left(k_i,k_{i+1}\right)}\left(t,y_{i+1}-y_i\right)},
\end{eqnarray*} where
$H_{\left(k_i,k_{i+1}\right)}\left(t,y_{i+1}-y_i\right)$ can be
calculated similarly to Eq.(\ref{13}).

\section{Random motion with reflecting boundaries}

\noindent Consider a set of real numbers  $\{y_i;
i=1,\dots,n\}\subset(0,b)$, where $b>0$ and $y_1<y_2<\dots<y_n$.
Let $S_1\left(t\right),S_2\left(t\right),\dots,S_n\left(t\right)$
be independent telegraph processes. It is assumed that all
processes have absolute velocity $v$ and parameter of switching
process $\lambda$ and starting form $y_i$ the process $S_i(t)$ has
equal probabilities of initial directions of the motion. We
suppose that $0$ and $b$ are two reflecting boundaries such that
if a process reaches boundary $0$ or $b$ then it changes velocity
direction to the opposite. Consider the family of particles with
trajectories $x(y_1,t),x(y_2,t),\dots,x(y_n,t)$, where every
$x(y_i,t)$ coincides respectively with processes $S_i(t)$ before
particle $i$ has first hard collision with another particle or
equivalently to the first intersection of the process $S_i(t)$
with another process. After the first hard collision of the
particle $x\left(y_i,t\right)$ with another particle, say
$x\left(y_j,t\right)$ they will switch the telegraph processes
that describe their trajectories so, $S_i\left(t\right)$ will
coincide with the trajectory of $x\left(y_j,t\right)$ and so on.

It is easily seen that the trajectories of the particles
$x\left(y_k,t\right)$, $k=1,2,\dots $ coincides with the order
statistics of $S_i(t)$, $i=1,2,\dots$,  as follows

\begin{equation} \label{15}
x\left(y_1,t\right)= S_{(1)}\left(t\right),
x\left(y_2,t\right)=S_{(2)}\left(t\right), \dots,
x\left(y_n,t\right)=S_{(n)}\left(t\right).
\end{equation}

\noindent Let us introduce the following distribution functions
$F_{y_r}\left(x\right)=P\left\{x\left(y_r,t\right)<x\right\}$.
Denote by $M_k^{(l)}$, $l=1,2,\dots,C_n^k$, different $k$ elements
subsets of the set $M=\{1,2,\dots,n\}$. \noindent It follows from
Eqs.(\ref{15}) that

\[F_{y_r}(x)=P\left\{x\left(y_r,t\right)<x\right\}=\sum_{k=r}^n\sum_{l=1}^{C_n^k}\prod_{i\in M_k^{(l)}}P(S_i
(t)<x)\prod_{j\in M\backslash M_k^{(l)}}P(S_j(t)\geq x).\]

\noindent For some particular cases we have

\[F_{y_1}\left(x\right)=P\left(x\left(y_1,t\right)<x\right)=1-\prod^n_{i=1}{P\left(S_i\left(t\right)\ge x\right)},\]

\begin{eqnarray*}
F_{y_{n-1}}\left(x\right)=P\left(x\left(y_{n-1},t\right)<x\right)=\sum^n_{k=1}\prod^n_{i=1,
i\ne k}
{P\left(S_i\left(t\right)<x\right)}P\left(S_k\left(t\right)\ge
x\right)\\
+\prod^n_{i=1}{P\left(S_i\left(t\right)<x\right)},
\end{eqnarray*}

\[F_{y_n}\left(x\right)=P\left(x\left(y_n,t\right)<x\right)=\prod^n_{i=1}{P\left(S_{i}\left(t\right)<x\right)}.\]

Let us study the limiting distribution of $S_k\left(t\right)$,
$k=1,\dots ,n$ as $t\to +\infty $. Denote by $N\left(t\right)$ the
number of Poisson events that have occurred in the interval
$\left(0,t\right)$ and let $s_j$, $j\ge 0$ be instants at which
Poisson events occur, and $s_0=0$. We assume that instants $s_j$
denote times of change of direction of $S_k\left(t\right)$.

\begin{lem} \label{l5}
Suppose that $f\left(x\right)$ is an integrable function on
$\left[0,b\right]$. Then
\[P\left({\mathop{\lim }_{T\to +\infty } \frac{1}{T}\int^T_0{f\left(S_k\left(t\right)\right)dt}\ }=\frac{1}{b}\int^b_0{f\left(x\right)dx}\right)=1.\]
\end{lem}

\noindent \textit{Proof} In the sequel, we will use the well-known
strong law of the large numbers for Poisson process

\begin{equation} \label{16}
P\left(\mathop{{\rm lim}}_{T\to +\infty
}\frac{N\left(T\right)}{T}=\lambda \right)=1.
\end{equation}

\noindent Since during the time $s_{j+1}-s_j$ the particle covers
the distance of $(s_{j+1}-s_j)v$ the number
$\left[\frac{(s_{j+1}-s_j )v}{2b}\right]$ is equal to the double
number of passages of segment $[0,b]$ by the particle.

Hence,

\begin{eqnarray}\label{17}
\nonumber {\mathop{\lim }_{T\to +\infty }
\frac{1}{T}\int^T_0{f\left(S_k\left(t\right)\right)dt}\
}&=&\mathop{{\rm lim}}_{T\to +\infty
}\frac{1}{T}\sum^{N\left(T\right)}_{i=0}{\int^{s_{j+1}}_{s_j}{f\left(S_k\left(t\right)\right)dt}}\\
&=&\mathop{{\rm lim}}_{T\to +\infty
}\frac{1}{T}\sum^{N\left(T\right)}_{i=0}{\left(\left[\frac{\left(s_{j+1}-s_j\right)v}{2b}\right]\frac{2}{v}\int^b_0f\left(x\right)dx+r_i\right)}
\ \ a.s.
\end{eqnarray}

\noindent where $r_i=\int_{u_i}^{u_i+\vartheta_i}f(S_k (t))dt$,
$u_i$, $\vartheta_i$ are independent random variables and $u_i$ is
uniformly distributed on $[0,2b]$, $\vartheta_i$ has the following
pdf
\begin{equation*}
g(t)=\frac{\lambda}{v}e^{-\frac{\lambda
t}{v}}\left(1-e^{-\frac{2\lambda b}{v}}\right)^{-1}I_{\{0\leq t
\leq 2b\}}.
\end{equation*}

\noindent Therefore,
\begin{eqnarray}\label{18}
\nonumber Er_i=E\int_{u_{i}}^{u_{i}+\vartheta_{i}}f(S_{k}(t))dt& =
&\frac{\lambda}{2bv \left(1-e^{-\frac{2\lambda
b}{v}}\right)}\int_{0}^{2b}dx\int_{0}^{2b}dp\int_{x}^{x+p}dt
f(S_{k}(t))e^{-\frac{\lambda p}{v}}\\
\nonumber & = & -\frac{1}{2b}\int_{0}^{2b}dx \int_{x}^{x+2b}dt
f(S_{k}(t))\frac{e^{-\frac{2\lambda
b}{v}}}{\left(1-e^{-\frac{2\lambda b}{v}}\right)}\\
& + & \frac{1}{2b\left(1-e^{-\frac{2\lambda
b}{v}}\right)}\int_{0}^{2b}dx\int_{0}^{2b}d p
f\left(S_{k}(x+p)\right)e^{-\frac{\lambda p}{v}}\\
\nonumber & = &-\frac{2e^{-\frac{2\lambda
b}{v}}}{\left(1-e^{-\frac{2\lambda
b}{v}}\right)}\int_{0}^{b}f(x)dx+\frac{1}{\lambda
b}\int_{0}^{b}f(x)dx.
\end{eqnarray}

\noindent The strong law of the large numbers for $\{r_i,i\geq
1\}$ implies
\begin{equation*}
\lim_{N\to +\infty }\frac{1}{N}\sum^{N}_{i=1}{r_i}=Er_i \ \ a.s.
\end{equation*}
Since ${\theta }_j=s_{j+1}-s_j$, $j=1,2,\dots $ are independent
exponentially distributed random variables, we have the following
strong law of the large numbers

\begin{eqnarray}\label{19}
\nonumber \lim_{N \to
+\infty}\frac{1}{N}\sum^{N}_{j=1}\left[\frac{(s_j-s_{j-1}
)v}{2b}\right]&=&E\left[\frac{(s_j-s_{j-1})v}{2b}\right]\\
&=&\sum^{\infty}_{n=1}n\left(e^{-\frac{2n\lambda
b}{v}}-e^{-\frac{2(n+1)\lambda
b}{v}}\right)=\frac{e^{-\frac{2\lambda
b}{v}}}{1-e^{-\frac{2\lambda b}{v}}}.
\end{eqnarray}
Combining Eqs.(\ref{16})-(\ref{19}), we get

\begin{eqnarray*}{\mathop{\lim }_{T\to +\infty }
\frac{1}{T}\int^T_0{f\left(S_k\left(t\right)\right)dt}}
=\mathop{\lim }_{T\to +\infty
}\frac{N\left(t\right)}{T}\frac{1}{N\left(T\right)}\sum^{N\left(T\right)}_{j=1}{\left(\int_{s_{j-1}}^{s_j}f(S_k(t))\right)}\\
= \frac{1}{b}\int_{0}^{b}f(x)dx \ \ \ a.s.
\end{eqnarray*}

\noindent This concludes the proof.

Therefore, the limiting distribution of $S_{k}\left(t\right)$ as
$t\to +\infty $ for all $k=1,\dots ,n$ is the uniform distribution
on $\left[0,b\right].$

\begin{lem} \label{l6}
Suppose that the initial distribution of a telegraph particle with
reflecting boundaries $0$ and b is uniform on $[0,b]$. Then it
remains uniform for all $t>0$.
\end{lem}

\noindent \textit{Proof} Denote by
$p\left(t,x\mathrel{\left|\vphantom{t,x
y_k}\right.\kern-\nulldelimiterspace}y_k\right)$ the probability
density of the process $S_k\left(t\right)$ position at time $t$.
It was shown in \cite{p15} that for $x\in \left[0,b\right]$

\[p\left(t,x\mathrel{\left|\vphantom{t,x y_k}\right.\kern-\nulldelimiterspace}y_k\right)=\frac{1}{b}+\frac{2}{b}e^{-\lambda t}\sum^{\infty }_{n=1}{\left\{\left[{\cosh  \left({\theta }_nt\right)\ }+\frac{\lambda }{{\theta }_n}{\sinh  {\theta }_nt\ }\right]{\cos  \left(\frac{\pi ny_k}{b}\right)\ }{\cos  \left(\frac{\pi nx}{b}\right)\ }\right\}},\]

\noindent where

\[{\theta }_n={\left({\lambda }^2-\frac{{\pi }^2v^2}{b^2}n^2\right)}^{{1}/{2}}.\]

\noindent It is easily seen that for any $t>0$ and $x\in
\left[0,b\right]$

\[p\left(t,x\right)=\frac{1}{b}\int^b_0{p\left(t,x\mathrel{\left|\vphantom{t,x y_k}\right.\kern-\nulldelimiterspace}y_k\right)}dy_k=\frac{1}{b}.\]

Now let us consider the system of processes
$\overline{S}_k\left(t\right)$ with the limiting distribution of
the respective processes $S_k\left(t\right)$, $k=1,2,\dots,n$.
According to Lemmas \ref{l5} and \ref{l6}, for each $t\ge 0$
processes $\overline{S}_k\left(t\right)$, $k=1,2,\dots n$ are
independent and having the uniform distribution on
$\left[0,b\right].$

For this case denote by $x_k\left(t\right)$, $k=1,2,\dots,n$
particles positions at time $t\geq 0$. It is easy to see that for
every $t\geq 0$ processes $x_k\left(t\right)$ are the order
statistics of $S_k\left(t\right)$, $k=1,2,\dots,n$, namely

\begin{eqnarray*}
x_1\left(t\right)=\overline{S}_{(1)}\left(t\right),
x_2\left(t\right)=\overline{S}_{(2)}\left(t\right),\dots,
x_n\left(t\right)=\overline{S}_{(n)}\left(t\right).
\end{eqnarray*}

Consider the following function

\[p\left(x\right)=P(S_k (t)<x)=\left\{\begin{array}{c}
\frac{x}{b},\ \ \ \ x\in \left[0,b\right],\  \\
\ 0,\ \ \ x\notin \left[0,b\right]. \end{array} \right.\]

\noindent It is easily verified that the distributions ${\pi
}_k\left(\cdot \right)$ of the particles positions
$x_k\left(t\right)$, $k\in \left\{1,2,\dots ,n\right\}$ are as
follows

\[{{\pi }_k\left(x\right)= P\left(x_k\left(t\right)<x\right)\ }=I_{p\left(x\right)}\left(k,n-k+1\right),\]
where

\begin{equation*}
I_{p\left(x\right)}\left(k,n-k+1\right)=\frac{\int^{p\left(x\right)}_0{t^{k-1}{\left(1-t\right)}^{n-k}dt}}{\int^1_0{t^{k-1}{\left(1-t\right)}^{n-k}dt}}.
\end{equation*}

Let us study the number of collisions $C_{\left(1,2,\dots
,n\right)}\left(0,t\right)$ of particles $x_k\left(t\right)$,
$k=1,2,\dots ,n$ during time $\left(0,t\right)$. It is easy to see
that $C_{\left(1,2,\dots ,n\right)}\left(0,t\right)$ is a number
of intersections of $\overline{S}_k\left(t\right)$, $k=1,2,\dots,
n$ for each $t>0$.

Denote by $I_{\left(k,l\right)}\left(0,t\right)$, $k\ne l$ the
number of intersection of processes $\overline{S}_k\left(t\right)$
and $\overline{S}_l\left(t\right)$ during time $\left(0,t\right)$.
Then it is easily verified that

\[C_{\left(1,2,\dots ,n\right)}\left(0,t\right)=\sum_{1\le k<l\le n}{I_{\left(k,l\right)}\left(0,t\right)}.\]
Therefore, let us analyze the distribution of
$I_{\left(k,l\right)}\left(0,t\right)$, $k\ne l$.

\noindent Since $\overline{S}_k\left(t\right)$ and
$\overline{S}_l\left(t\right)$ have the uniform distribution on
$\left[0,b\right]$ the probability of their intersections
$I_{\left(k,l\right)}\left(t,t+\triangle t\right)$ during
$\left(t,t+\triangle t\right)$ satisfies the following
inequalities for $x\in \left(a,b\right)$

\begin{eqnarray} \label{20}
\nonumber
\frac{1}{4}P\left(\left|\overline{S}_k\left(t\right)-\overline{S}_l\left(t\right)\right|\le
2\triangle tv\right)e^{-2\lambda \triangle t} \le
P\left(N_{\left(k,l\right)}\left(t,t+\triangle t\right)\ge
1\right)\\
\le
P\left(\left|\overline{S}_k\left(t\right)-\overline{S}_l\left(t\right)\right|
\le 2\triangle tv\right).
\end{eqnarray}
By using
$P\left(\left|\overline{S}_k\left(t\right)-\overline{S}_l\left(t\right)\right|\le
2\triangle tv\right)=O\left(\triangle t\right)$ and
$\frac{1}{4}P\left(\left|\overline{S}_k\left(t\right)-\overline{S}_l\left(t\right)\right|\le
2\triangle tv\right)e^{-2\lambda \triangle t}=O\left(\triangle
t\right)$, we get

\begin{equation} \label{21}
P\left(I_{\left(k,l\right)}\left(t,t+\triangle t\right)\ge
1\right)=O\left(\triangle t\right).
\end{equation}
It is easily verified that for $n\ge 2$

\begin{eqnarray} \label{22}
\nonumber P\left(I_{\left(k,l\right)}\left(t,t+\triangle
t\right)=n\right)\le
P\left(\left|\overline{S}_k\left(t\right)-\overline{S}_l\left(t\right)\right|\le
2\triangle tv\right){\left(1-e^{-\lambda \triangle
t}\right)}^{2\left(n-1\right)}\\
+ \nonumber
P\left(\{\overline{S}_k\left(t\right),\overline{S}_l\left(t\right)\in
\left[0,2\triangle tv\right]\}\cup
\{\overline{S}_k\left(t\right),\overline{S}_l\left(t\right)\in
\left[b-2\triangle tv,b\right]\}\right){\left(1-e^{-\lambda
\triangle
t}\right)}^{\left(n-1\right)}\\
\le \frac{4\triangle
tv}{b}{\left(\lambda \triangle
t\right)}^{2\left(n-1\right)}+2\frac{{\left(2\triangle
tv\right)}^2}{b^2}{\left(\lambda \triangle
t\right)}^{\left(n-1\right)}.
\end{eqnarray}

\noindent Taking into account Eqs.(\ref{21}),(\ref{22}), we
conclude that

\[P\left(I_{\left(k,l\right)}\left(t,t+\triangle t\right)=1\right)=O\left(\triangle t\right).\]
Therefore, for $x\in \left(a,b\right)$

\begin{eqnarray} \label{23}
\nonumber O\left(\triangle
t\right)&=&\frac{1}{4}P\left(\left|\overline{S}_k\left(t\right)-\overline{S}_l\left(t\right)\right|\le
2\triangle tv\right)e^{-2\lambda \triangle t}
\\&\le&
\nonumber EI_{\left(k,l\right)}\left(t,t+\triangle t\right)\le
P\left(\left|\overline{S}_k\left(t\right)-\overline{S}_l\left(t\right)\right|\le
2\triangle tv\right)\\
\nonumber &+&\frac{4\triangle tv}{b}\sum_{n\ge
1}{n\left({\left(\lambda \triangle t\right)}^{2\left(n-1\right)}
+\frac{2\triangle tv}{b}{\left(\lambda \triangle
t\right)}^{\left(n-1\right)}\right)}
\\&=&\frac{4\triangle
tv}{b}+o\left(\triangle t\right).
\end{eqnarray}
It is easily seen the additive property of
$EI_{\left(k,l\right)}\left(t_1,t_2\right)$: for any $s\in
\left(t_1,t_2\right)$

\[EI_{\left(k,l\right)}\left(t_1,t_2\right)=EI_{\left(k,l\right)}\left(t_1,s\right)+EI_{\left(k,l\right)}\left(s,t_2\right).\]
Hence, there exists a constant $c>0$ such that

\[EI_{\left(k,l\right)}\left(0,t\right)=ct.\]
This implies that

\[EC_{\left(1,2,\dots ,n\right)}\left(0,t\right)=\frac{n\left(n-1\right)}{2}ct.\]

It follows from (\ref{20}) and (\ref{23}) the following estimation
for the factor $c$

\[\frac{v}{b}\le c\le \frac{4v}{b}.\]

\section{acknowledgments}

The author expresses deep gratitude to Professor A.A. Dorogovtsev
for graciously dedicating so much of his time, patience, and
knowledge in guiding me through the writing of this article.

\noindent
\noindent \eject

\end{document}